\documentclass[10pt]{amsart}
\usepackage{amsmath,amssymb,latexsym,hyperref,url,graphicx,stmaryrd}

\newcommand{\Z}{\mathbb{Z}}
\newcommand{\freq}{\textnormal{freq}}
\newcommand{\weight}{\textnormal{weight}}
\newcommand{\minratio}{\operatorname{minratio}}
\newcommand{\one}{\overset{\bigcurlyvee}{1\vphantom{\displaystyle{\sum}}}}
\newcommand{\two}{\multicolumn{2}{c}{\overset{\underbrace{}}{2\vphantom{\displaystyle{\sum}}}}}

\newtheorem*{lemma}{Lemma}
\newtheorem*{theorem}{Theorem}

\begin{document}

\title{Bounds on the frequency of $1$\\in the Kolakoski word}
\author{Elizabeth J. Kupin}
\author{Eric S. Rowland}
\address{
	Department of Mathematics\\
	Rutgers University\\
	Piscataway, NJ 08854, USA
}
\date{September 26, 2008}

\begin{abstract}
We use a method of Goulden and Jackson to bound $\freq_1(K)$, the limiting frequency of $1$ in the Kolakoski word $K$.  We prove that $|\freq_1(K) - 1/2| \leq 17/762$, assuming the limit exists, and establish the semi-rigorous bound $|\freq_1(K) - 1/2| \leq 1/46$.
\end{abstract}

\maketitle
\markboth{E. Kupin and E. Rowland}{Bounds on the frequency of $1$ in the Kolakoski word}

\section{Introduction}\label{introduction}

The Kolakoski word is an infinite sequence of $1$'s and $2$'s that is equal to its own run length sequence:

\[
\begin{array}{cccccccccccccccccccccc}
K = & 2 & 2 & 1 & 1 & 2 & 1 & 2 & 2 & 1 & 2 & 2 & 1 & 1 & 2 & 1 & 1 & 2 & 2 & 1 & 2 & \cdots \\
K = & \two & \two & \one & \one & \two & \one & \two & \two & \one & \two & \two & \one & \one & \cdots
\end{array}
\]
Up to the choice of the first term, $K$ is defined uniquely by this property.  Beginning with $1$ instead of $2$ produces the word $1K$, which was introduced by Kolakoski \cite{kolakoski,kolakoski-ucoluk}.

Let $o_n$ be the number of $1$'s occurring in the first $n$ terms of $K$, and let
$$
	\freq_1(K) := \lim_{n \to \infty} \frac{o_n}{n}.
$$
It was conjectured by Dekking \cite{dekking} that this limit exists and equals $1/2$.  Kimberling's web page \cite{kimberling}, where this conjecture is listed among several others, is responsible for its popularity.  In this paper we use the Goulden--Jackson cluster method to give bounds on $\freq_1(K)$ consistent with the conjecture.  In particular, we prove the following.

\begin{theorem}
If the limiting frequency $\freq_1(K)$ of $1$ in the Kolakoski word exists, then
\[
	\left|\freq_1(K) - \frac{1}{2}\right| \leq \frac{17}{762} \approx 0.0223097 .
\]
\end{theorem}

This method was explored in a different setting by Chv\'atal \cite{chvatal}, who produced this bound and several better bounds, reducing the difference from $1/2$ to
\[
	\left|\freq_1(K) - \frac{1}{2}\right| \leq \frac{35}{41754} \approx 0.0008382.
\]
We thank Jean-Paul Allouche for pointing out Chv\'atal's paper to us.

\section{The Goulden--Jackson cluster method}\label{gouldenjackson}

\subsection{Description}\label{description}

The Goulden--Jackson cluster method \cite{goulden-jackson} is an efficient way of counting the number of words $w$ on a given alphabet such that no subword of $w$ appears in a given set $S$.  We say that $w$ \emph{avoids} $S$.

Here we use an extension of the method by Noonan and Zeilberger \cite{noonan-zeilberger} that tracks the frequency of the letters in a word.  Define the \emph{weight} of a word $w$ to be
\[
	\weight(w) = x_1^{|w|_1} \, x_2^{|w|_2} \, t_{\vphantom{1}}^{|w|},
\]
where $|w|_\alpha$ is the number of occurrences of $\alpha$ in $w$ and $|w|$ is the length of $w$.  Let $W$ be the set of words on $\{1, 2\}$ that avoid $S$.  Let the weight of $W$ be
\[
	\weight(W) = \sum_{w \in W} \weight(w) = \sum_{n=0}^\infty p_n(x_1,x_2) t^n,
\]
where $p_n(x_1,x_2)$ is a polynomial in $x_1$ and $x_2$ that carries information about the set of length-$n$ words avoiding $S$.  The Goulden--Jackson algorithm computes $\weight(W)$ as a rational expression in $x_1$, $x_2$, and $t$.  We refer the reader to the papers cited above for details of the algorithm.

\subsection{Avoided subwords}\label{avoidedsubwords}

To use the Goulden--Jackson method we must find words that never appear as subwords of $K$.  We accomplish this by capitalizing on the fact that if $w$ is a subword of $K$ then the run length sequence of $w$ is also a subword of $K$.  This means that if $w$ is not a subword of $K$ then any word whose run length sequence contains $w$ is also not a subword of $K$. We start by observing that the word $3$ does not occur in $K$ because $K$ is a word on $\{1, 2\}$.  Therefore, no word with $3$ in its run length sequence can be a subword of $K$ either; in particular, $111$ and $222$ cannot be subwords of $K$.  

Now that we know that $K$ avoids $111$ and $222$, we know that no word with $111$ or $222$ in its run length sequence can occur in $K$.  Namely, $K$ avoids $12121$ and $21212$ (since their run length sequences contain $111$), and $K$ also avoids $112211$ and $221122$ (since their run length sequences are $222$).

There is a subtlety here, which is that $111$ is the run length sequence of the words $212$ and $121$, yet these words both do appear in $K$. However, they only occur as part of the larger words $22122$ and $11211$, and these word have run length sequences $212$, not $111$. We pad $212$ with $1$'s on both ends to ensure that the run length sequence contains $111$, and similarly we pad $121$ with $2$'s. This padding is necessary whenever the run length sequence begins or ends with $1$.

\begin{figure}
	\includegraphics{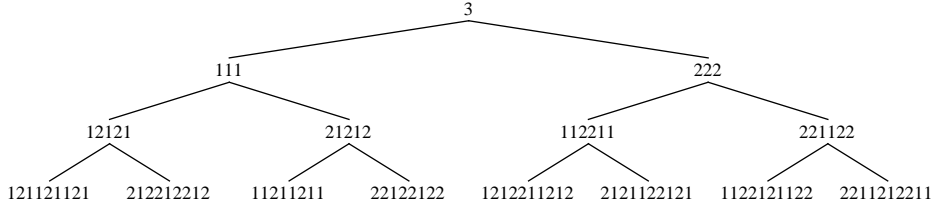}
	\caption{The first four generations in an infinite tree of words that the Kolakoski word avoids.}
	\label{wordtree}
\end{figure}

We iterate this process to obtain additional words that $K$ avoids, producing the tree in Figure~\ref{wordtree}.  Define $S_d$ be the set of words in the tree in levels $1$ through $d$ (i.e., not including the root, $3$). There are $2^{d+1}-2$ words in $S_d$.

This approach to producing words avoided by $K$ is symmetric with respect to interchanging $1$ and $2$, so if we use all words in $S_d$ it follows that $p_n(x_1,x_2)$ is symmetric in $x_1$ and $x_2$.  Because of this symmetry, all our bounds have the form $|\freq_1(K) - 1/2| \leq \epsilon$.  Experiments with asymmetric word sets have not improved upon the bounds obtained with symmetric sets, so we do not pursue them here.

\section{Results}\label{results}

We have three different (but closely related) ways of using the Goulden--Jackson method to produce bounds on $\freq_1(K)$. When we can compute the full generating function $\weight(W)$, the denominator gives a bound directly. For large sets $S$ computing the generating function as a rational expression is not computationally feasible; in this case we resort to computing the first few terms of the series expansion.  Each term of the series provides a bound on $\freq_1(K)$, although in general these bounds are not as good as the ones we get from the denominator. Finally, by examining many terms of the series we can often experimentally determine a closed form for the bounds being produced, which, after taking a limit, gives an improved bound that is semi-rigorous.

\subsection{Bounds from the denominator}\label{denominator}

From the term $p_n(x_1,x_2) t^n$ we can determine the minimum number of $1$'s that occur in an $n$-letter word avoiding $S$; this is the minimum degree in $x_1$ of this polynomial.  Let
\[
	\minratio \, \sum_{i=1}^k c_i x_1^{o_i} x_2^{n_i-o_i} t^{n_i} = \min_{1 \leq i \leq k} \frac{o_i}{n_i}
\]
for $c_i \in \Z$ and $n_i \geq 1$.

If $\weight(W) = N/(1-D)$ for some polynomials $N$ and $D$, then $\weight(W) = \sum_{n=0}^\infty N D^n$, and
\[
	\minratio N D^n \to \minratio D
\]
as $n \to \infty$.  Thus the denominator of $\weight(W)$ dictates the asymptotic behavior of $\minratio p_n(x_1,x_2) t^n$.

For example, using the set $S_1=\{111, 222\}$ produces the generating function
\[
	\weight(W) = \frac{\left(x_1^2 t^2 + x_1 t + 1\right) \left(x_2^2t^2 +x_2 t + 1\right)}{1 - x_1^2 x_2^2 t^4 - x_1^2 x_2 t^3 - x_1x_2^2 t^3 - x_1 x_2 t^2}.
\]
Here $\minratio D = 1/3$, and the maximum ratio is $2/3$.  Therefore if the limit exists we have $|\freq_1(K) - 1/2| \leq 1/6$.

The $\minratio$ for $S_2$ is also $1/3$ despite additional words.  However, using $S_3$ produces the denominator
\begin{multline*}
	1 + x_1^{18} x_2^{18} t^{36} - x_1^{16} x_2^{17} t^{33} - x_1^{17} x_2^{16} t^{33} - x_1^{15} x_2^{15} t^{30} + 3 x_1^{12} x_2^{12} t^{24} \\
	+ x_1^{10} x_2^{11} t^{21} + x_1^{11} x_2^{10} t^{21} + x_1^8 x_2^{10} t^{18} + x_1^9 x_2^9 t^{18} + x_1^{10} x_2^8 t^{18} - x_1^7 x_2^8 t^{15} \\
	- x_1^8 x_2^7 t^{15} - 2 x_1^6 x_2^6 t^{12} - x_1^5 x_2^5 t^{10} - 2 x_1^4 x_2^5 t^9 - 2 x_1^5 x_2^4 t^9 - x_1^4 x_2^4 t^8
\end{multline*}
with $\minratio D = 4/9$, giving $|\freq_1(K) - 1/2| \leq 1/18$.

\subsection{Bounds from series terms}\label{seriesterms}

Computing $\weight(W)$ as a rational expression requires solving a system of linear equations, and this system is large when there are many words in $S$.  Therefore, to compute improved bounds we use a modified algorithm, available in the function \texttt{wGJseries} in Zeilberger's package \texttt{DAVID\_IAN} \cite{davidian}, that computes only the first $N$ terms of the series.  The following proposition says that each term puts a bound on $\freq_1(K)$.  The idea is that bounding the number of $1$'s in every length-$n$ subword of $K$ produces a bound that extends to all of $K$.  Recall that $o_m$ is the number of $1$'s in the first $m$ terms of $K$.

\begin{lemma}
If $a \leq |w|_1 \leq b$ for every length-$n$ subword $w$ of $K$, then $a/n \leq \freq_1(K) \leq b/n$ if the limit exists.
\end{lemma}
\begin{proof}
For an arbitrary $m$, we have $m = q n + r$ for some $q \in \Z$ and $0 \leq r < n$.  Partition the first $m$ terms of $K$ into $q$ consecutive blocks of length $n$, leaving a remainder block of length $r$. The number of $1$'s in each block of length $n$ is at most $b$, so $o_m \leq q b + r$.  Similarly, $o_m \geq q a$. Therefore for all $m$ we have
$$
	\frac{qa}{m} \leq  \frac{o_m}{m} \leq \frac{qb+r}{m}.
$$
Substituting $q = \frac{m-r}{n} $ and letting $m \to \infty$ we get that
$$
	\frac{a}{n} \leq \lim_{m \to \infty} \frac{o_m}{m} \leq \frac{b}{n}
$$
if the limit in question exists.
\end{proof}

For example, we compute $\weight(W)$ with $S_1 = \{111, 222\}$ out to term $N=5$ to be

\begin{multline*}
	1 + \left(x_1 + x_2\right) t +\left(x_1^2 + 2x_1x_2 + x_2^2\right) t^2 + \left(3x_1^2x_2 + 3x_1x_2^2\right) t^3 \\
	+ \left(2x_1^3x_2 + 6x_1^2x_2^2 + 2x_1x_2^3\right) t^4 + \left(x_1^4x_2 + 7x_1^3x_2^2 + 7x_1^2x_2^3 + x_1x_2^4\right) t^5.
\end{multline*}
From the coefficient of $t^3$, we conclude that $1 \leq |w|_1 \leq 2$ for every word of length $3$ avoiding $111$ and $222$.  This information gives the bound $|\freq_1(K)-1/2| \leq 1/6$, which in this case is the same bound obtained from the denominator.  While we could have used any coefficient in the series expansion to get a bound, the bounds we get from the coefficients of $t^4$ and $t^5$ are actually worse.

Performing similar computations on $S_d$ for larger $d$ produces better bounds.  The following table gives the best bound $\epsilon(n) = 1/2 - \minratio p_n(x_1,x_2) t^n$ achieved among the first $N$ terms.  Computing $N = 800$ terms for $d=5$ took a day and a half.
\[
\begin{array}{ccccc}
d & |S_d| & N & n & \epsilon(n) \\ \hline
1 & 2 & 200 & 3 & 1/6 \\
2 & 6 & 200 & 3 & 1/6 \\
3 & 14 & 200 & 9 & 1/18 \\
4 & 30 & 500 & 498 & 17/498 \\
5 & 62 & 800 & 762 & 17/762 \\
6 & 126 & 600 & 555 & 5/222
\end{array}
\]
The best bound here is $|\freq_1(K) - 1/2| \leq 17/762$, provided by $d = 5$.

For $1 \leq d \leq 3$ the bounds achieved are best possible for these word sets; indeed they are the same bounds obtained from $\minratio D$ for $\weight(W)$.  For $d \geq 4$, computing more terms will produce increasingly better bounds, although for a fixed $d$ the bounds approach $1/2 - \minratio D$, as discussed in the following section.  Likewise, using more words should produce better bounds, although this increases the computation time.

\subsection{Implied bounds}\label{semi-rigorous}

In fact the sequence of minimum degrees of $x_1$ in $p_n(x_1,x_2)$ has the simple structure of a linear quasi-polynomial for sufficiently large $n$.  Moreover, the successive maxima of the sequence $\minratio p_n(x_1,x_2) t^n$ eventually occur in just one of the residue classes.

Having computed several terms for $d=1$ it is not difficult to guess that for $n \geq 1$ the minimum degree of $x_1$ in $p_n(x_1,x_2)$ is given by the linear quasi-polynomial
\[
	\begin{cases}
		\frac{n}{3}	& \text{if $n \equiv 0 \bmod 3$,} \\
		\frac{n - 1}{3}	& \text{if $n \equiv 1 \bmod 3$,} \\
		\frac{n - 2}{3}	& \text{if $n \equiv 2 \bmod 3$.}
	\end{cases}
\]
Therefore $\minratio p_n(x_1,x_2) t^n \to 1/3$, and in fact the limit is attained every three terms beginning at $n=3$.  The sequence of minimum degrees for $d=2$ is identical to that for $d=1$.

For $d=3$, $n \minratio p_n(x_1,x_2) t^n$ is given by
\[
	\begin{cases}
		4 \cdot \frac{n}{9}		& \text{if $n \equiv 0 \bmod 9$,} \\
		4 \cdot \frac{n - 1}{9}		& \text{if $n \equiv 1 \bmod 9$,} \\
		4 \cdot \frac{n - 2}{9}		& \text{if $n \equiv 2 \bmod 9$,} \\
		4 \cdot \frac{n - 3}{9} + 1	& \text{if $n \equiv 3 \bmod 9$,} \\
		4 \cdot \frac{n - 4}{9} + 1	& \text{if $n \equiv 4 \bmod 9$,} \\
		4 \cdot \frac{n - 5}{9} + 1	& \text{if $n \equiv 5 \bmod 9$,} \\
		4 \cdot \frac{n - 6}{9} + 2	& \text{if $n \equiv 6 \bmod 9$,} \\
		4 \cdot \frac{n - 7}{9} + 2	& \text{if $n \equiv 7 \bmod 9$,} \\
		4 \cdot \frac{n - 8}{9} + 3	& \text{if $n \equiv 8 \bmod 9$.}
	\end{cases}
\]
The limit, $4/9$, is first attained at $n=9$, producing $\epsilon = 1/18$.

For higher values of $d$, the limit is not attained by any term.  The sequence of minimum ratios for $d=4$ and $d=5$ are eventually linear quasi-polynomials.  Using $d=6$ iterations of words, one finds that at least the first $600$ terms in the series have the same $\minratio$ as those for $d=5$, with the exception of $n=62$; therefore the same eventual quasi-polynomial seems to hold.

For $d=4$ the quasi-polynomial has modulus $15$.  For residue class $n \equiv i \bmod 15$ the main term is $7 \cdot \frac{n-i}{15}$, and the constant terms for $i = 0, 1, \dots, 14$ are
\[
	-1, -1, 0, 1, 1, 1, 2, 2, 2, 3, 4, 4, 5, 5, 5.
\]
The successive maxima are $(7 m + 1)/(15 m + 3)$ for $m \geq 2$ and occur at terms $n = 15 m + 3$.  Thus the limit is $7/15$, producing $\epsilon = 1/30$.

For $d=5$ the modulus is $69$.  The main term is $33 \cdot \frac{n-i}{69}$ for $n \equiv i \bmod 69$, and the constant terms are
\begin{multline*}
	-1, -1, 0, 1, 1, 1, 2, 2, 2, 3, 4, 4, 5, 5, 5, 6, 6, 7, 8, 8, 8, 9, 9, 9, 10, 11, 11, \\
	12, 12, 13, 13, 14, 14, 15, 15, 15, 16, 17, 17, 18, 18, 18, 19, 19, 20, 21, 21, 21, \\
	22, 22, 22, 23, 24, 24, 25, 25, 25, 26, 26, 27, 28, 28, 28, 29, 29, 30, 31, 31, 31.
\end{multline*}
The successive maxima are $(33 m + 1)/(69 m + 3)$ for $m \geq 3$ and occur at terms $n = 69 m + 3$; for example, $m = 11$ produces the best rigorous bound
\[
	\epsilon = \frac{1}{2} - \frac{33 \cdot 11 + 1}{69 \cdot 11 + 3} = \frac{17}{762}.
\]
Therefore most probably $\minratio D = 33/69$ for $\weight(W)$ in this case, and
\[
	\left|\freq_1(K) - \frac{1}{2}\right| \leq \frac{1}{46} \approx 0.0217391 .
\]

\end{document}